\documentclass[11pt,leqno]{article}

\usepackage[active]{srcltx}

\usepackage{amsfonts,latexsym,amsmath,amssymb,amsthm}
\usepackage{fullpage}
\newtheorem{theorem}{Theorem}[section]
\newtheorem{lemma}[theorem]{Lemma}

\newtheorem{proposition}[theorem]{Proposition}

\newtheorem{remark}[theorem]{Remark}
\newtheorem{question}[theorem]{Question}

\theoremstyle{definition}
\newtheorem{definition}[theorem]{Definition}

\theoremstyle{remark}

\newtheorem*{note*}{Note}

\numberwithin{equation}{section}

\makeatletter
\newcommand{\rank}{\mathop{\operator@font rank}}
\newcommand{\conv}{\mathop{\operator@font conv}}
\newcommand{\vol}{\mathop{\operator@font vol}}
\newcommand{\onetagright}{\tagsleft@false}
\makeatother

\newcommand{\ls}{\leqslant}
\newcommand{\gr}{\geqslant}
\renewcommand{\epsilon}{\varepsilon}
\newcommand{\prend}{$\quad \hfill \Box$}
\usepackage{color}

\usepackage{hyperref}

\begin{document}
\small

\title{\bf Inequalities for the surface area of projections of convex bodies}

\medskip

\author{Apostolos Giannopoulos, Alexander Koldobsky and Petros Valettas}

\date{}

\maketitle

\begin{abstract}
\footnotesize We provide general inequalities that compare the surface area $S(K)$ of a convex body $K$ in ${\mathbb R}^n$
to the minimal, average or maximal surface area of its hyperplane or lower dimensional projections. We discuss the
same questions for all the quermassintegrals of $K$. We examine separately the dependence of the constants
on the dimension in the case where $K$ is in some of the classical positions or $K$ is a projection body.
Our results are in the spirit of the hyperplane problem, with sections replaced by projections and volume by
surface area.
\end{abstract}

\section{Introduction}

The starting point of this article are two inequalities of the second named author about the surface area of
hyperplane projections of projection bodies. In \cite{Koldobsky-2013} it was proved that if $Z$ is a projection body in ${\mathbb R}^n$ then
\begin{equation}\label{eq:intro-1}|Z|^{\frac{1}{n}}\,\min_{\xi\in S^{n-1}}S(P_{\xi^{\perp }}(Z))\ls b_nS(Z),\end{equation}
where $S(A)$ denotes the surface area of $A$ and
\begin{equation}\label{eq:bn-def}b_n=\frac{(n-1)\omega_{n-1}}{n\omega_n^{\frac{n-1}{n}}},\end{equation}
where $\omega_m$ is the volume of the Euclidean unit ball $B_2^m$ in ${\mathbb R}^m$. Note that $b_n\simeq 1$
and that \eqref{eq:intro-1} is sharp; there is equality if $Z=B_2^n$. Conversely, in \cite{Koldobsky-2015}
it was proved that if $Z$ is a projection body in ${\mathbb R}^n$ which is a dilate of a body in isotropic position, then
\begin{equation}\label{eq:intro-2}|Z|^{\frac{1}{n}}\,\max_{\xi\in S^{n-1}}S(P_{\xi^{\perp }}(Z))\gr c (\log n)^{-2}S(Z),\end{equation}
where $c>0$ is an absolute constant.

Our first aim is to discuss similar inequalities for the surface area of hyperplane projections of an arbitrary convex body $K$ in ${\mathbb R}^n$.
In what follows, we denote by $\partial_K$ the minimal surface area parameter of $K$, defined by
\begin{equation}\label{eq:partial}\partial_K:=\min\Big\{ S(T(K))/|T(K)|^{\frac{n-1}{n}}:T\in GL(n)\Big\}.\end{equation}
It is known that $c\sqrt{n}\ls \partial_K\ls c^{\prime }n$ for every convex body $K$ in ${\mathbb R}^n$, where $c,c^{\prime }>0$
are absolute constants (see Section 2 for definitions, references and background information).

Our analogue of \eqref{eq:intro-1} is the following theorem.

\begin{theorem}\label{th:intro-1}There exists an absolute constant $c_1>0$ such that, for every convex body $K$ in ${\mathbb R}^n$,
\begin{equation}\label{eq:intro-3}|K|^{\frac{1}{n}}\,\min_{\xi\in S^{n-1}}S(P_{\xi^{\perp }}(K))\ls \frac{2b_n\partial_K}{n\omega_n^{\frac{1}{n}}}\,S(K)\ls\frac{c_1\partial_K}{\sqrt{n}}\,S(K).\end{equation}
\end{theorem}

This inequality is sharp e.g. for the Euclidean unit ball. Note that $c_1\partial_K/\sqrt{n}\ls c\sqrt{n}$ for every convex
body $K$ in ${\mathbb R}^n$, and hence we have the general upper bound
\begin{equation}|K|^{\frac{1}{n}}\,\min_{\xi\in S^{n-1}}S(P_{\xi^{\perp }}(K))\ls c\sqrt{n}\,S(K).\end{equation}
Our method employs an estimate for the minimal volume of a hyperplane projection of $K$: one has
\begin{equation}\label{eq:intro-4}\min_{\xi\in S^{n-1}}|P_{\xi^{\perp }}(K)|\ls c\sqrt{n}\,|K|^{\frac{n-1}{n}}\end{equation}
for an absolute constant $c>0$.

Assuming that $K$ is in the minimal surface area position we have a converse of Theorem \ref{th:intro-1}:

\begin{theorem}\label{th:intro-lower-min}Let $K$ be a convex body in ${\mathbb R}^n$ which is in the minimal surface area position. Then,
\begin{equation}\label{eq:intro-lower-min}|K|^{\frac{1}{n}}\,\min_{\xi\in S^{n-1}}S(P_{\xi^{\perp }}(K))\gr \frac{c}{\sqrt{n}}\,S(K),\end{equation}
where $c>0$ is an absolute constant.
\end{theorem}

The estimate of Theorem \ref{th:intro-lower-min} is sharp; we provide an example in which the two quantities in \eqref{eq:intro-lower-min}
are of the same order, using extremal (with respect to minimal hyperplane projections) bodies of minimal surface area that were constructed
in \cite{Markessinis-Paouris-Saroglou-2012}.

In the case where $K$ is a projection body, one can see that \eqref{eq:intro-4} holds true
with $c\sqrt{n}$ replaced by $b_n$ (see Section 3). This leads to an alternative proof of \eqref{eq:intro-1} with a weaker (by a factor of $2$) constant.

\begin{theorem}\label{th:intro-2}Let $Z$ be a projection body in ${\mathbb R}^n$. Then,
\begin{equation}\label{eq:intro-5}|Z|^{\frac{1}{n}}\,\min_{\xi\in S^{n-1}}S(P_{\xi^{\perp }}(Z))\ls 2b_n\,S(Z).\end{equation}
\end{theorem}

It should be noted that there are convex bodies which are not projection bodies but their minimal surface area parameter $\partial_K$ is
of the order of $\sqrt{n}$; an example is given by $B_1^n$, the unit ball of $\ell_1^n$. On the other hand, there exist projection bodies
whose minimal surface area parameter is of the order of $n$; an example is given by the cube. Thus, the estimates of Theorem \ref{th:intro-1}
and Theorem \ref{th:intro-2} complement each other. In our next result we replace $\min\,S(P_{\xi^{\perp }}(K))$ by the expectation
of $S(P_{\xi^{\perp }}(K))$ on the sphere.

\begin{theorem}\label{th:intro-3}Let $K$ be a convex body in ${\mathbb R}^n$. Then,
\begin{equation}\label{eq:intro-6}|K|\,\int_{S^{n-1}}S(P_{\xi^{\perp }}(K))\,d\sigma (\xi )\ls \frac{2(n-1)\omega_{n-1}}{n^2\omega_n}S(K)^2
\ls \frac{c_2}{\sqrt{n}}S(K)^2,\end{equation}
where $c_2>0$ is an absolute constant.
\end{theorem}

A consequence of Theorem \ref{th:intro-3} is that if $K$ is in some of the {\it classical positions} (minimal surface area, isotropic or John's position,
or it is symmetric and in L\"{o}wner's position) then
\begin{equation}\label{eq:intro-7}|K|^{\frac{1}{n}}\,\int_{S^{n-1}}S(P_{\xi^{\perp }}(K))\,d\sigma (\xi )\ls c\sqrt{n}\,S(K).\end{equation}
The reason is that, in all these cases, the surface area of $K$ satisfies an inequality of the form $S(K)\ls cn|K|^{\frac{n-1}{n}}$ (see Section 2
for a brief description of the classical positions of a convex body and for a proof of this last assertion).

\smallskip

Passing to lower bounds, our analogue of \eqref{eq:intro-2} is the following theorem.

\begin{theorem}\label{th:intro-4}Let $K$ be a convex body in ${\mathbb R}^n$. Then,
\begin{equation}\label{eq:intro-8}\int_{S^{n-1}}S(P_{\xi^{\perp }}(K))\,d\sigma (\xi )\gr c_3\,S(K)^{\frac{n-2}{n-1}},\end{equation}
where $c_3>0$ is an absolute constant.
\end{theorem}

A consequence of Theorem \ref{th:intro-4} is that if $K$ is in the minimal surface area, minimal mean width, isotropic, John or L\"{o}wner position, then
\begin{equation}\label{eq:intro-9}|K|^{\frac{1}{n}}\,\int_{S^{n-1}}S(P_{\xi^{\perp }}(K))\,d\sigma (\xi )\gr c\,S(K),\end{equation}
where $c>0$ is an absolute constant. In particular,
\begin{equation}\label{eq:intro-10}|K|^{\frac{1}{n}}\,\max_{\xi\in S^{n-1}}S(P_{\xi^{\perp }}(K))\gr c\,S(K).\end{equation}
Note that \eqref{eq:intro-10} is stronger than \eqref{eq:intro-2}; moreover, for bounds of this type there is no need to assume that $K$ is a projection body.
In fact, our proof of Theorem \ref{th:intro-4} shows that \eqref{eq:intro-9} continues to hold as long as
\begin{equation}\label{eq:intro-11}S(K)^{\frac{1}{n-1}}\ls c|K|^{\frac{1}{n}}\end{equation}
for an absolute constant $c>0$. This is a mild condition which is satisfied not only by the classical positions but also
by all {\it reasonable} positions of $K$.

All these inequalities are proved in Section 4. Our main tools are a result from \cite{Giannopoulos-Hartzoulaki-Paouris-2002} stating that
\begin{equation}\label{eq:intro-12}\frac{S(P_{\xi^{\perp }}(K))}{|P_{\xi^{\perp }}(K)|}\ls \frac{2(n-1)}{n}\frac{S(K)}{|K|}\end{equation}
for every convex body $K$ in ${\mathbb R}^n$ and any $\xi\in S^{n-1}$, estimates from \cite{Giannopoulos-Papadimitrakis-1999}
for the volume of the projection body of a convex body in terms of its minimal surface area parameter, and Aleksandrov's inequalities.
For the class of projection bodies, we prove and use the following sharp estimate (Lemma \ref{lem:zon-1}): if $Z$ is a projection body in ${\mathbb R}^n$ then
\begin{equation}\label{eq:intro-13}\min_{\xi\in S^{n-1}}|P_{\xi^{\perp }}(Z)|\ls \frac{nb_n}{n-1}\,|Z|^{\frac{n-1}{n}}.\end{equation}

In Section 5 we study the same questions for the quermassintegrals $V_{n-k}(K)=V((K, n-k), (B_2^n,k))$ of a convex body $K$ and the
corresponding quermassintegrals of its hyperplane projections. We obtain the following estimates:
\begin{enumerate}
\item[(i)] For every $1\ls p\ls n-2$ we have
\begin{equation}\label{eq:intro-14}|K|^{\frac{1}{n}}\,\min_{\xi\in S^{n-1}}V_{n-1-p}(P_{\xi^{\perp }}(K))\ls
\frac{(p+1)\omega_{n-1}\partial_K}{n\omega_n}\,V_{n-p}(K)\ls \frac{c(p+1)\partial_K}{\sqrt{n}}\,V_{n-p}(K)\end{equation}
and
\begin{equation}\label{eq:intro-15}|K|^{\frac{1}{n}}\,\int_{S^{n-1}}V_{n-1-p}(P_{\xi^{\perp }}(K))\,d\sigma (\xi )\ls
\frac{(p+1)\omega_{n-1}}{n\omega_n}\frac{S(K)}{|K|^{\frac{n-1}{n}}}\,V_{n-p}(K).\end{equation}
\item[(ii)] If $Z$ is a projection body in ${\mathbb R}^n$ then, for every $1\ls p\ls n-2$ we have
\begin{equation}\label{eq:intro-16}|Z|^{\frac{1}{n}}\,\min_{\xi\in S^{n-1}}V_{n-1-p}(P_{\xi^{\perp }}(Z))\ls (p+1)b_n\,V_{n-p}(Z).\end{equation}
\item[(iii)] If $K$ is in the minimal surface area, isotropic or John's position, or it is symmetric and in L\"{o}wner's position, then, for every $1\ls p\ls n-2$ we have
\begin{equation}\label{eq:intro-17}|K|^{\frac{1}{n}}\,\int_{S^{n-1}}V_{n-1-p}(P_{\xi^{\perp }}(K))\,d\sigma (\xi )\ls
c(p+1)\sqrt{n}\,V_{n-p}(K).\end{equation}
\item[(iv)] For every $1\ls p\ls n-2$ we have
\begin{equation}\int_{S^{n-1}}V_{n-1-p}(P_{\xi^{\perp }}(K))\,d\sigma (\xi )\gr
\frac{\omega_{n-1}}{\omega_n^{\frac{n-1-p}{n-p}}}\,[V_{n-p}(K)]^{\frac{n-1-p}{n-p}}.\end{equation}
\item[(v)] If $K$ is in the minimal surface area, isotropic or John's position,
or it is symmetric and in L\"{o}wner's position then, for every $1\ls p\ls n-2$ we have
\begin{equation}|K|^{\frac{1}{n}}\int_{S^{n-1}}V_{n-1-p}(P_{\xi^{\perp }}(K))\,d\sigma (\xi )\gr
\frac{\omega_{n-1}c_0^{\frac{p}{n-p}}}{\omega_n^{\frac{n-1-p}{n-p}}}\,V_{n-p}(K)\gr \left (\frac{c_1}{n}\right )^{\frac{p}{2(n-p)}}\,V_{n-p}(K).\end{equation}
\end{enumerate}
The proofs employ the same tools as in the surface area case. The main additional ingredient is a generalization of \eqref{eq:intro-12} to
subspaces of arbitrary dimension and quermassintegrals of any order, proved in \cite{Fradelizi-Giannopoulos-Meyer-2003}: If $K$ is a
convex body in ${\mathbb R}^n$ and $0\ls p\ls k\ls n$, then, for every $F\in G_{n,k}$,
\begin{equation}\label{eq:intro-18}\frac{V_{n-p}(K)}{|K|}\gr\frac{1}{\binom{n-k+p}{n-k}}\frac{V_{k-p}(P_F(K))}{|P_F(K)|}.\end{equation}
This inequality allows us to obtain further generalizations of the results of Section 4; we can compare the surface area of a convex body $K$
to the minimal, average or maximal surface area of its lower dimensional projections $P_F(K)$, $F\in G_{n,k}$, for any given $1\ls k\ls n-1$.
This is done in Section 6.

There are several questions that arise from this work and we hope that the reader might find them interesting; these are stated explicitely
throughout the text.

\section{Notation and background}

We work in ${\mathbb R}^n$, which is equipped with a Euclidean structure $\langle\cdot ,\cdot\rangle $. We denote by $\|\cdot \|_2$
the corresponding Euclidean norm, and write $B_2^n$ for the Euclidean unit ball and
$S^{n-1}$ for the unit sphere. We denote the unit ball of $\ell_p^n$ by $B_p^n$, $1\ls p\ls\infty $. In particular, we also write $Q_n$
for the cube $B_{\infty }^n=[-1,1]^n$ and $C_n=\left [-\frac{1}{2},\frac{1}{2}\right ]^n$ for the cube of volume $1$.
Volume is denoted by $|\cdot |$. We write $\omega_n$ for the volume of $B_2^n$ and $\sigma $ for the rotationally invariant probability
measure on $S^{n-1}$. The Grassmann manifold $G_{n,k}$ of $k$-dimensional subspaces of ${\mathbb R}^n$ is equipped with the Haar probability
measure $\nu_{n,k}$. For every $1\ls k\ls n-1$ and $F\in G_{n,k}$ we write $P_F$ for the orthogonal projection from $\mathbb R^{n}$ onto $F$,
and we set $B_F=B_2^n\cap F$ and $S_F=S^{n-1}\cap F$. Finally, we write $\overline{A}$ for the homothetic image of volume 1 of a symmetric
convex body $A\subseteq \mathbb R^n$, i.e. $\overline{A}:=|A|^{-\frac{1}{n}}A$.

The letters $c,c^{\prime }, c_1, c_2$ etc. denote absolute positive constants which may change from line to line. Whenever we write
$a\simeq b$, we mean that there exist absolute constants $c_1,c_2>0$ such that $c_1a\ls b\ls c_2a$.  Also, if $K,L\subseteq \mathbb R^n$
we will write $K\simeq L$ if there exist absolute constants $c_1, c_2>0$ such that $ c_{1}K\subseteq L \subseteq c_{2}K$.

We refer to the books \cite{Gardner-book} and \cite{Schneider-book} for basic facts from the Brunn-Minkowski theory and to the book
\cite{AGA-book} for basic facts from asymptotic convex geometry. We also refer to \cite{BGVV-book} for more information on isotropic convex bodies.

\medskip

\noindent {\bf 2.1. Convex bodies}. A convex body in ${\mathbb R}^n$ is a compact convex subset $K$ of ${\mathbb R}^n$ with non-empty interior. We say that $K$ is
symmetric if $x\in K$ implies that $-x\in K$, and that $K$ is centered if its barycenter $\frac{1}{|K|}\int_Kx\,dx $ is at the origin. The support function of a
convex body $K$ is defined by $h_K(y)=\max \{\langle x,y\rangle :x\in K\}$, and the mean width of $K$ is
\begin{equation}\label{eq:not-1}w(K)=\int_{S^{n-1}}h_K(\theta )\,d\sigma (\theta ). \end{equation}
The circumradius of $K$ is the quantity $R(K)=\max\{ \|x\|_2:x\in K\}$ i.e. the smallest $R>0$ for which $K\subseteq RB_2^n$.
If $0\in {\rm int}(K)$ then we write $r(K)$ for the inradius of $K$ (the largest $r>0$
for which $rB_2^n\subseteq K$) and we define the polar body $K^{\circ }$ of $K$ by
\begin{equation}\label{eq:not-2}K^{\circ }:=\{ y\in {\mathbb R}^n: \langle x,y\rangle \ls 1 \;\hbox{for all}\; x\in K\}. \end{equation}
The volume radius of $K$ is the quantity ${\rm vrad}(K)=\left (|K|/|B_2^n|\right )^{1/n}$.
Integration in polar coordinates shows that if the origin is an interior point of $K$ then the volume radius of $K$ can be expressed as
\begin{equation}\label{eq:not-3}{\rm vrad}(K)=\left (\int_{S^{n-1}}\|\theta\|_K^{-n}\,d\sigma (\theta )\right)^{1/n},\end{equation}
where $\|\theta \|_K=\min\{ t>0:\theta \in tK\}$. We also define
\begin{equation}M(K)=\int_{S^{n-1}}\|\theta\|_K\,d\sigma (\theta ).\end{equation}

\medskip

\noindent {\bf 2.2. Mixed volumes}. From Minkowski's fundamental theorem we know that if $K_1,\ldots ,K_m$ are non-empty, compact convex
subsets of ${\mathbb R}^n$, then the volume of $t_1K_1+\cdots +t_mK_m$ is a homogeneous polynomial of degree $n$ in
$t_i>0$. That is,
\begin{equation}\label{eq:not-4}|t_1K_1+\cdots +t_mK_m|=\sum_{1\ls i_1,\ldots ,i_n\ls m}
V(K_{i_1},\ldots ,K_{i_n})t_{i_1}\cdots t_{i_n},\end{equation}
where the coefficients $V(K_{i_1},\ldots ,K_{i_n})$ are chosen to be invariant under permutations of their arguments. The coefficient $V(K_1,\ldots ,K_n)$ is the mixed
volume of $K_1,\ldots ,K_n$. In particular, if $K$ and $C$ are two convex bodies in ${\mathbb R}^n$
then the function $|K+tC|$ is a polynomial in $t\in [0,\infty )$:
\begin{equation}\label{eq:not-5}|K+tC|=\sum_{j=0}^n \binom{n}{j} V_{n-j}(K,C)\;t^j,\end{equation}
where $V_{n-j}(K,C)= V((K,n-j),(C,j))$ is the $j$-th mixed volume of $K$ and $C$ (we use  the notation $(C,j)$ for $C,\ldots ,C$ $j$-times).
If $C=B_2^n$ then we set $V_{n-j}(K):=V_{n-j}(K,B_2^n)=V((K, n-j), (B_2^n, j))$; this is the $j$-th quermassintegral of $K$.
Note that
\begin{equation}\label{eq:not-6}V_{n-1}(K,C)={\frac{1}{n}} \lim_{t\to 0^+}{\frac{|K+tC|-|K|}{t}},\end{equation}
and by the Brunn-Minkowski inequality we see that
\begin{equation}\label{eq:not-7}V_{n-1}(K,C)\gr |K|^{\frac{n-1}{n}}|C|^{1/n}\end{equation}
for all $K$ and $C$ (this is Minkowski's first inequality). The mixed volume $V_{n-1}(K,C)$ can be expressed as
\begin{equation}\label{eq:not-8}V_{n-1}(K,C)={\frac{1}{n}}\int_{S^{n-1}}h_C(\theta )d\sigma_K(\theta ),\end{equation}
where $\sigma_K$ is the surface area measure of $K$. In particular, the surface area of $K$ satisfies
\begin{equation}\label{eq:not-9}S(K)=nV_{n-1}(K).\end{equation}
We will also use the Aleksandrov inequalities: if $K$ is a convex body in ${\mathbb R}^n$ then the sequence
\begin{equation}\label{eq:aleksandrov-1}Q_k(K)=\left
(\frac{1}{\omega_k}\int_{G_{n,k}}|P_F(K)|\,d\nu_{n,k}(F)\right )^{\frac{1}{k}}\end{equation}is decreasing in $k$.
This is a consequence of the Aleksandrov-Fenchel inequality (see \cite{Burago-Zalgaller-book} and \cite{Schneider-book}).
In particular, for every $1\ls k\ls n-1$ we have
\begin{equation}\label{eq:aleksandrov-2}\left (\frac{|K|}{\omega_n}\right )^{\frac{1}{n}}\ls \left (\frac{1}{\omega_k}\int_{G_{n,k}}|P_F(K)|\,d\nu_{n,k}(F)\right )^{\frac{1}{k}}\ls w(K).\end{equation}

\medskip

\noindent {\bf 2.3. Classical positions}. Let $K$ be a centered convex body in ${\mathbb R}^n$.
We introduce the classical positions of $K$ that we are going to discuss; we set the notation and provide some background information.

\smallskip

\noindent {\it Minimal surface area position}. We say that $K$ has minimal surface area if
$S(K)\ls S(T(K))$ for every $T\in SL(n)$. Petty
(\cite{Petty-1961}, see also \cite{Giannopoulos-Papadimitrakis-1999}) proved that $K$ has minimal surface
area if and only if the measure $\sigma_K$ satisfies the isotropic condition
\begin{equation}\label{eq:not-10}S(K)=n\int_{S^{n-1}}\langle \xi ,\theta\rangle^2d\sigma_K (\theta )\end{equation}
for every $\xi\in S^{n-1}$. From the isoperimetric inequality we know that $S(K)\gr n\omega_n^{\frac{1}{n}}|K|^{\frac{n-1}{n}}$.
The reverse isoperimetric inequality of K.~Ball \cite{Ball-reverse} implies that if $K$ has
minimal surface area and volume $1$ then $S(K)\ls S(C_n)=2n$ in the symmetric case and $S(K)\ls
S(\Delta_n)\ls c_0n$ in the not necessarily symmetric case, where $\Delta_n$ is a regular simplex
of volume $1$ in ${\mathbb R}^n$ and $c_0>0$ is an absolute constant.

\smallskip

\noindent {\it Minimal mean width position}. We say that $K$ is in minimal mean width position if $w(K)\ls
w(T(K))$ for every $T\in SL(n)$. It was proved in \cite{Giannopoulos-VMilman-2000} that
$K$ has minimal mean width if and only if
\begin{equation}\label{eq:not-11}w(K)=n\int_{S^{n-1}}\langle \xi ,\theta\rangle^2h_K(\theta )d\sigma (\theta )\end{equation}
for every $\xi\in S^{n-1}$. From results of Figiel-Tomczak, Lewis and Pisier (see \cite[Chapter 6]{AGA-book})
it follows that if a symmetric convex body $K$ in ${\mathbb R}^n$ has minimal mean width then
\begin{equation}\label{eq:not-12}M(K)w(K)\ls c_1\log (d_K+1)\end{equation}
where $d_K:=d(K,B_2^n)$ is the Banach-Mazur distance from $K$ to $B_2^n$ and $c_1>0$ is an absolute constant.
If we assume that $|K|=1$ then $w(K)\ls c\sqrt{n}\log (d_K+1)$. Then, a simple argument shows that any convex
body of volume $1$ in ${\mathbb R}^n$ that has minimal mean width satisfies a similar bound: $w(K)\ls c\sqrt{n}\log n$.

\smallskip

\noindent {\it Isotropic position}. For every centered convex body $K$ of volume $1$ in ${\mathbb R}^n$
and any $q\gr 1$ we define
\begin{equation}\label{eq:not-13}I_q(K)=\left (\int_K\|x\|_2^qdx\right)^{1/q}.\end{equation}
We say that $K$ is in the isotropic position if $I_2(K)\ls I_2(T(K))$ for every $T\in SL(n)$. This is equivalent to the existence of a
constant $L_K>0$ such that
\begin{equation}\label{eq:not-14}\int_K\langle x,\xi\rangle^2dx =L_K^2\end{equation}
for every $\xi\in S^{n-1}$. It is known that if $K$ is centered then
\begin{equation}\label{eq:not-15}\int_K\langle x,\xi\rangle^2dx\simeq |K\cap\xi^{\perp }|^{-2}\end{equation}
for every $\xi\in S^{n-1}$. Therefore, if $K$ is isotropic we see that all hyperplane sections $K\cap\xi^{\perp }$
of $K$ have volume equal (up to an absolute constant) to $L_K^{-1}$.

\smallskip

\noindent {\it John and L\"{o}wner position.} We say that a convex body $K$ is in John's position if the
ellipsoid of maximal volume inscribed in $K$ is a multiple of the Euclidean unit ball $B_2^n$. We say that a convex body $K$ is in L\"{o}wner's
position if the ellipsoid of minimal volume containing $K$ is a multiple of the Euclidean unit ball
$B_2^n$. One can check that this holds true if and only if $K^{\circ }$ is in John's position.
The volume ratio of a centered convex body $K$ in ${\mathbb R}^n$ is the quantity
\begin{equation}\label{eq:not-16}{\rm vr}(K)=\inf\left\{ \left(\frac{|K|}{|{\cal E}|}\right )^{\frac{1}{n}}:{\cal E}\;\hbox{is an ellipsoid and}\;
{\cal E}\subseteq K\right\}.\end{equation}
The outer volume ratio of a centered convex body $K$ in ${\mathbb R}^n$ is the quantity ${\rm ovr}(K)={\rm vr}(K^{\circ })$.
K.~Ball proved in \cite{Ball-reverse} that if $K$ is in John's position then ${\rm vr}(K)\ls {\rm vr}(C_n)\simeq \sqrt{n}$ in the symmetric case and
${\rm vr}(K)\ls {\rm vr}(\Delta_n)\simeq \sqrt{n}$ in the not necessarily symmetric case; in fact, the reverse
isoperimetric inequality follows from this fact.

\medskip

\noindent {\bf 2.4. Surface area and inradius}. Let $K$ be a centered convex body in ${\mathbb R}^n$.
Recall that the inradius $r(K)$ of $K$ is the largest $r>0$ for which $rB_2^n\subseteq K$. Using the monotonicity of mixed volumes we may write
\begin{equation}\label{eq:not-17}S(K) = nV_{n-1}(K,B_2^n)\ls n V_{n-1}\Big(K,\frac{1}{r(K)} K\Big).\end{equation}
Since the mixed volumes are homogeneous with respect to each of their arguments and $V(K,\ldots ,K)=|K|$,
we have the following general estimate for the surface area $S(K)$ of $K$.

\begin{lemma}\label{lem:surface-inradius}
Let $K$ be a convex body in $\mathbb R^n$ with $0\in {\rm int}(K)$. Then,
\begin{equation}\label{eq:not-18}S(K)\ls\frac{n|K|}{r(K)}.\end{equation}
\end{lemma}

Using Lemma \ref{lem:surface-inradius} we obtain upper bounds for the surface area of a body which
is in isotropic, John's or L\"{o}wner's position.

\begin{proposition}\label{prop:surface-IJL}Let $K$ be a centered convex body of volume $1$ in ${\mathbb R}^n$.
\begin{enumerate}
\item[{\rm (i)}] If $K$ is isotropic then $S(K) \ls cn/L_K\ls c^{\prime }n$, where $c, c^{\prime }>0$ are absolute constants.
\item[{\rm (ii)}] If $K$ is in minimal surface area position or in John's position
then $S(K)\ls cn$, where $c>0$ is an absolute constant.
\item[{\rm (iii)}] If $K$ is symmetric and in L\"{o}wner's position
then $S(K)\ls cn$, where $c>0$ is an absolute constant.
\item[{\rm (iv)}] If $K$ is symmetric and in the minimal mean width position
then $S(K)\ls cn\log n$, where $c>0$ is an absolute constant.
\end{enumerate}
\end{proposition}

\noindent {\it Proof.} The inclusion $L_KB_2^n\subseteq K$ for an isotropic symmetric convex body $K$ in ${\mathbb R}^n$
is clear since \begin{equation*}h_{K}(u)=\|\langle\cdot
,u\rangle\|_{L^{\infty }(K)}\gr\|\langle\cdot
,u\rangle\|_{L^2(K)}=L_K\end{equation*} for every $u\in S^{n-1}$. This shows that
$r(K)\gr L_K$ in this case. If $K$ is centered but not necessarily symmetric, then we still have
$h_K(u)\gr cL_K$: to see this, we use the fact that $c\max\{ |K\cap (t\theta +\theta^{\perp })|:t\in {\mathbb R}\}\ls |K\cap\theta^{\perp }|$
(see \cite[Chapter 2]{BGVV-book}) and then write
\begin{align*}L_K^{-1}h_K(u) &\gr c_1h_K(u)|K\cap\theta^{\perp }|\gr c_2\int_0^{\infty }|K\cap (t\theta +\theta^{\perp })|\,dt\\
&=c_2|\{x\in K:\langle x,\theta\rangle\gr 0\}|\gr c_3,\end{align*}where $c_3>0$ is an absolute constant (the last inequality follows
from Gr\"{u}nbaum's lemma, see \cite[Chapter 2]{BGVV-book}). To conclude the proof we recall that $L_K\gr c$ for any convex body $K$
in ${\mathbb R}^n$.

Assume that $K$ is in John's position. Then, using the volume ratio estimate we see that
\begin{equation}\frac{\sqrt{n}}{r(K)}\simeq \left (\frac{|K|}{|r(K)B_2^n|}\right )^{\frac{1}{n}}= {\rm vr}(K)\ls c\sqrt{n},\end{equation}
which implies that $r(K)\gr c$, and hence $S(K)\ls c^{-1}n$. It follows that if $K$ is in minimal surface area position we also have $S(K)\ls
c^{-1}n$.

Next, assume that $K$ is symmetric and in L\"{o}wner's position; this time we use the fact
that $R(K)\ls \sqrt{n}r(K)$ by John's theorem, and then
\begin{equation}1=|K|^{1/n}\ls |R(K)B_2^n|^{1/n}\ls cR(K)/\sqrt{n}\ls cr(K).\end{equation}
Finally, if $K$ is symmetric and in the minimal mean width position we can use the direct estimate
\begin{equation}R(K^{\circ})\ls c\sqrt{n}w(K^{\circ })=c\sqrt{n}M(K)\ls c^{\prime }\log n\end{equation}
which is a consequence of \eqref{eq:not-12} and of the fact that $w(K)\gr c\sqrt{n}$ by Urysohn's
inequality. This shows that $r(K)=1/R(K^{\circ })\gr c/\log n$, and (iv) follows. \prend

\medskip

\noindent {\it Note.} The example of the cube $C_n$ shows that the bounds (i), (ii) and (iii) of Proposition \ref{prop:surface-IJL}
are sharp up to an absolute constant.

\section{Projections of projection bodies}

A zonoid is the limit of Minkowski sums of line segments in the Hausdorff metric. Equivalently, a
symmetric convex body $Z$ is a zonoid if and only if its polar body is the unit ball of an $n$-dimensional subspace of an $L_1$-space;
i.e. if there exists a positive measure $\mu $ (the supporting measure of $Z$) on $S^{n-1}$ such that
\begin{equation}\label{eq:zon-1}h_Z(x)=\| x\|_{Z^{\circ }}=\frac{1}{2}\int_{S^{n-1}}|\langle x,y\rangle |d\mu (y).\end{equation}
The class of origin-symmetric zonoids coincides with the class of projection bodies. Recall that the projection body $\Pi K$ of a convex body $K$ is the
symmetric convex body whose support function is defined by
\begin{equation}\label{eq:zon-2}h_{\Pi K} (\xi )=|P_{\xi^{\perp } }(K)|, \qquad \xi\in S^{n-1}.\end{equation}
From Cauchy's formula
\begin{equation}\label{eq:zon-3}|P_{\xi^{\perp } }(K)|={\frac{1}{2}}\int_{S^{n-1}}|\langle \xi ,\theta \rangle |\;d\sigma_K(\theta ),\end{equation}
where $\sigma_K$ is the surface area measure of $K$, we see that the projection body of $K$ is a zonoid whose
supporting measure is $\sigma_K$. Minkowski's existence theorem implies that, conversely, every origin-symmetric zonoid
is the projection body of some symmetric convex body in ${\mathbb R}^n$.
Moreover, if we denote by ${\cal C}_n$ the class of origin-symmetric convex bodies
and by ${\cal Z}_n$ the class of origin-symmetric zonoids, Aleksandrov's uniqueness theorem shows that the Minkowski map
$\Pi : {\cal C}_n \to {\cal Z}_n$ with $K\mapsto \Pi K$, is
injective. Note also that ${\cal Z}_n$ is invariant under invertible
linear transformations and closed in the Hausdorff metric.

Let $K$ be a convex body of volume $1$ in ${\mathbb R}^n$. Then, the volume of $\Pi K$ and of its polar body $\Pi^{\ast }K$ satisfy
the bounds (see \cite{Giannopoulos-Papadimitrakis-1999})
\begin{equation}\label{eq:zon-4}\left (\frac{\partial_K}{n}\right )^n\ls |\Pi K|
\ls\omega_n\left (\frac{\omega_{n-1}\partial_K}{n\omega_n}\right )^n\end{equation}
and
\begin{equation}\label{eq:zon-5}\omega_n\left (\frac{n\omega_n}{\omega_{n-1}\partial_K}\right )^n\ls |\Pi^{\ast }K|
\ls \frac{4^nn^n}{n!\partial_K^n}.\end{equation}
All these inequalities are sharp as one can see from the examples of the ball and the cube.

\smallskip

Our next lemma provides an estimate for the volume of the minimal hyperplane projection of a zonoid.

\begin{lemma}\label{lem:zon-1}Let $Z$ be a zonoid in ${\mathbb R}^n$. Then,
\begin{equation}\label{eq:zon-6}\min_{\xi\in S^{n-1}}|P_{\xi^{\perp }}(Z)|\ls \frac{nb_n}{n-1}\,|Z|^{\frac{n-1}{n}}.\end{equation}
\end{lemma}

\noindent {\it Proof.} We write $Z=\Pi K$ for some convex body $K$. Recall the volume formula for zonoids
\begin{equation}\label{eq:zon-7}|Z|=V_{n-1}(Z,\Pi K)=V_{n-1}(K,\Pi Z)=\frac{1}{n}\int_{S^{n-1}}h_{\Pi Z}(\xi )\,d\sigma_K(\xi )=
\frac{1}{n}\int_{S^{n-1}}|P_{\xi^{\perp }}(Z)|\,d\sigma_K(\xi ).\end{equation}
Therefore,
\begin{equation}\label{eq:zon-8}|Z|\gr \frac{S(K)}{n}\min_{\xi \in S^{n-1}}|P_{\xi^{\perp }}(Z)|.\end{equation}
On the other hand, by \eqref{eq:zon-4} we have
\begin{equation}|Z|^{\frac{1}{n}}=|\Pi K|^{\frac{1}{n}}\ls\frac{\omega_{n-1}\partial_K}{n\omega_n^{\frac{n-1}{n}}}|K|^{\frac{n-1}{n}}\ls \frac{b_nS(K)}{n-1},\end{equation}
where we used the definition \eqref{eq:bn-def} of the constant $b_n$ and the fact that $S(K)\gr \partial_K|K|^{\frac{n-1}{n}}$ by the definition
\eqref{eq:partial} of the minimal surface area parameter $\partial_K$. Then, 
\begin{equation}\label{eq:zon-9}\frac{S(K)}{n}\min_{\xi \in S^{n-1}}|P_{\xi^{\perp }}(Z)|\ls |Z|^{\frac{1}{n}}|Z|^{\frac{n-1}{n}}\ls \frac{b_nS(K)}{n-1}|Z|^{\frac{1}{n}},\end{equation}
and the result follows. \prend

\medskip

Since every projection of a zonoid is a zonoid, a simple induction argument leads to the following general result.

\begin{theorem}\label{th:zon-2}Let $Z$ be a zonoid in ${\mathbb R}^n$. Then, for every $1\ls k\ls n-1$ we have
\begin{equation}\label{eq:zon-10}\min_{F\in G_{n,k}}|P_F(Z)|\ls \frac{nb_n^{n-k}}{k}\,|Z|^{\frac{k}{n}}.\end{equation}
\end{theorem}

\begin{question}\label{question-zon}\rm The example of the ball shows that Lemma \ref{lem:zon-1} is sharp. It would be interesting to
establish the precise constant in Theorem \ref{th:zon-2}.
\end{question}

\section{Surface area of hyperplane projections}

Our generalization of \eqref{eq:intro-1} is in terms of the minimal surface area parameter $\partial_K$ of $K$.

\begin{theorem}\label{th:hyper-1}Let $K$ be a convex body in ${\mathbb R}^n$. Then,
\begin{equation}\label{eq:hyper-1}|K|^{\frac{1}{n}}\,\min_{\xi\in S^{n-1}}S(P_{\xi^{\perp }}(K))\ls
\frac{2b_n\partial_K}{n\omega_n^{\frac{1}{n}}}\,S(K)\ls\frac{c_1\partial_K}{\sqrt{n}}\,S(K),\end{equation}
where $c_1>0$ is an absolute constant and $\partial_K$ is the minimal surface area parameter of $K$. Therefore, in general we have
that
\begin{equation}\label{eq:hyper-1-gen}|K|^{\frac{1}{n}}\,\min_{\xi\in S^{n-1}}S(P_{\xi^{\perp }}(K))
\ls c\sqrt{n}\,S(K).\end{equation}
\end{theorem}

The main ingredient in the proof is the next result from \cite{Giannopoulos-Hartzoulaki-Paouris-2002}.

\begin{lemma}\label{lem:GHP}If $K$ is a convex body in ${\mathbb R}^n$ then
\begin{equation}\label{eq:GHP}\frac{S(P_{\xi^{\perp }}(K))}{|P_{\xi^{\perp }}(K)|}\ls \frac{2(n-1)}{n}\frac{S(K)}{|K|}\end{equation}
for every $\xi\in S^{n-1}$.
\end{lemma}

\noindent {\bf Proof of Theorem \ref{th:hyper-1}.} From \eqref{eq:GHP} we have
\begin{equation}\label{eq:hyper-2}|K|\,S(P_{\xi^{\perp }}(K))\ls \frac{2(n-1)}{n}\,S(K)\,|P_{\xi^{\perp }}(K)|\end{equation}
for every $\xi\in S^{n-1}$.  Therefore,
\begin{equation}\label{eq:hyper-3}|K|\,\min_{\xi\in S^{n-1}}S(P_{\xi^{\perp }}(K))\ls \frac{2(n-1)}{n}\,S(K)\,\min_{\xi\in S^{n-1}}|P_{\xi^{\perp }}(K)|.\end{equation}
Next, we observe that
\begin{equation}\label{eq:hyper-4}\min_{\xi\in S^{n-1}}|P_{\xi^{\perp }}(K)|=\min_{\xi\in S^{n-1}}h_{\Pi K}(\xi )=r(\Pi K).\end{equation}
Since
\begin{equation}\label{eq:hyper-5}|\Pi K|^{\frac{1}{n}}\ls \frac{\omega_{n-1}\partial_K}{n\omega_n^{\frac{n-1}{n}}}\,|K|^{\frac{n-1}{n}}\end{equation}
we get
\begin{equation}\label{eq:rPiK}r(\Pi K)\ls {\rm vrad}(\Pi K)\ls \frac{\omega_{n-1}\partial_K}{n\omega_n}\,|K|^{\frac{n-1}{n}}.\end{equation}
Going back to \eqref{eq:hyper-3} we see that
\begin{equation}\label{eq:hyper-6}|K|\,\min_{\xi\in S^{n-1}}S(P_{\xi^{\perp }}(K))\ls
\frac{2(n-1)\omega_{n-1}\partial_K}{n^2\omega_n}\,S(K)\,|K|^{\frac{n-1}{n}},\end{equation}
and this proves \eqref{eq:hyper-1}. \prend

\medskip

\begin{question}\label{question-1}\rm It would be interesting to decide whether there exist convex bodies $K$ such that
\begin{equation}|K|^{\frac{1}{n}}\,\min_{\xi\in S^{n-1}}S(P_{\xi^{\perp }}(K))\gr c\sqrt{n}\,S(K).\end{equation}
This would show that Theorem \ref{th:hyper-1} is asymptotically sharp. Note that in the case of the Euclidean ball one has
\begin{equation}|B_2^n|^{\frac{1}{n}}\,\min_{\xi\in S^{n-1}}S(P_{\xi^{\perp }}(B_2^n))\gr c\,S(B_2^n).\end{equation}
\end{question}

We can prove an inequality which is reverse to \eqref{eq:hyper-1} for any convex body $K$ which is in the minimal surface area position, using the following fact
(see Theorem 1.2 in \cite{Giannopoulos-Hartzoulaki-Paouris-2002}): for any convex body $K$ in ${\mathbb R}^n$ and any
$\xi\in S^{n-1}$ one has
\begin{equation}\frac{V_{n-2}(K)}{2V_{n-1}(K)}\ls\frac{V_{n-2}(P_{\xi^{\perp }}(K))}{|P_{\xi^{\perp }}(K)|}.\end{equation}
Note that
\begin{equation}V_{n-2}(K)=\omega_n\,[Q_{n-2}(K)]^{n-2}\gr \omega_n\,[Q_{n-1}(K)]^{n-2}=\omega_n^{\frac{1}{n-1}}[V_{n-1}(K)]^{\frac{n-2}{n-1}},\end{equation}
while
\begin{equation}nV_{n-1}(K)=S(K), \quad (n-1)V_{n-2}(P_{\xi^{\perp }}(K))=S(P_{\xi^{\perp }}(K))\end{equation}
and if $K$ is in the minimal surface area position by \eqref{eq:not-10} we also have
\begin{equation}|P_{\xi^{\perp }}(K)|=\frac{1}{2}\int_{S^{n-1}}|\langle u,\xi\rangle |\,d\sigma_K(u)\gr
\frac{1}{2}\int_{S^{n-1}}\langle u,\xi\rangle^2\,d\sigma_K(u)=\frac{S(K)}{2n}=\frac{V_{n-1}(K)}{2}.\end{equation}
Combining the above we get
\begin{equation}S(P_{\xi^{\perp }}(K))\gr \frac{(n-1)\omega_n^{\frac{1}{n-1}}}{4}[V_{n-1}(K)]^{\frac{n-2}{n-1}}.\end{equation}
Therefore, we have
\begin{equation}|K|^{\frac{1}{n}}\,\min_{\xi\in S^{n-1}}S(P_{\xi^{\perp }}(K)) \gr \frac{(n-1)\omega_n^{\frac{1}{n-1}}}{4}[V_{n-1}(K)]^{\frac{n-2}{n-1}}\,|K|^{\frac{1}{n}}=
\frac{(n-1)\omega_n^{\frac{1}{n-1}}}{4n^{\frac{n-2}{n-1}}}[S(K)]^{\frac{n-2}{n-1}}\,|K|^{\frac{1}{n}}.\end{equation}
Since $S(K)=\partial_K\,|K|^{\frac{n-1}{n}}$, we get:
\begin{equation}|K|^{\frac{1}{n}}\,\min_{\xi\in S^{n-1}}S(P_{\xi^{\perp }}(K))\gr
\frac{(n-1)\omega_n^{\frac{1}{n-1}}}{4n^{\frac{n-2}{n-1}}\partial_K^{\frac{1}{n-1}}}\,S(K).\end{equation}
This proves the following:

\begin{theorem}\label{th:lower-min}Let $K$ be a convex body in ${\mathbb R}^n$ which is in the minimal surface area position. Then,
\begin{equation}|K|^{\frac{1}{n}}\,\min_{\xi\in S^{n-1}}S(P_{\xi^{\perp }}(K))\gr \frac{c}{\sqrt{n}}\,S(K),\end{equation}
where $c>0$ is an absolute constant.
\end{theorem}

\begin{remark}\rm Theorem \ref{th:lower-min} is sharp. In \cite{Markessinis-Paouris-Saroglou-2012} it is proved that
there exists an unconditional convex body $K_0$ of volume $1$ in
${\mathbb R}^n$ which has minimal surface area and satisfies
\begin{equation}\min_{\xi\in S^{n-1}}|P_{\xi^{\perp }}(K_0)|\ls\frac{c_1}{\sqrt{n}},\end{equation}
where $c_1>0$ is an absolute constant. From \eqref{eq:hyper-3} we see that
\begin{equation}|K_0|\,\min_{\xi\in S^{n-1}}S(P_{\xi^{\perp }}(K_0))\ls \frac{2(n-1)}{n}\,S(K_0)\,\min_{\xi\in S^{n-1}}|P_{\xi^{\perp }}(K_0)|
\ls \frac{2c_1}{\sqrt{n}}\,S(K_0),\end{equation}
and since $|K_0|=1$ we get:
\end{remark}

\begin{proposition}There exists an unconditional convex body $K_0$ of volume $1$ in
${\mathbb R}^n$ which has minimal surface area and satisfies
\begin{equation}|K_0|^{\frac{1}{n}}\,\min_{\xi\in S^{n-1}}S(P_{\xi^{\perp }}(K_0))\ls \frac{c}{\sqrt{n}}\,S(K_0),\end{equation}
where $c>0$ is an absolute constant.
\end{proposition}

Next, assume that $K=Z$ is a zonoid. Repeating the previous argument, but using now Lemma \ref{lem:zon-1} in order to estimate $r(\Pi Z)$, we obtain
an alternative proof of \eqref{eq:intro-1} (with a weaker, by a factor of $2$, constant).

\begin{theorem}\label{th:hyper-2}Let $Z$ be a zonoid in ${\mathbb R}^n$. Then,
\begin{equation}\label{eq:hyper-7}|Z|^{\frac{1}{n}}\,\min_{\xi\in S^{n-1}}S(P_{\xi^{\perp }}(Z))\ls
2b_n\,S(Z).\end{equation}
\end{theorem}

\noindent {\it Proof.} We have
\begin{equation}\label{eq:hyper-8}|Z|\,\min_{\xi\in S^{n-1}}S(P_{\xi^{\perp }}(Z)) \ls \frac{2(n-1)}{n}\,S(Z)\,\min_{\xi\in S^{n-1}}|P_{\xi^{\perp }}(Z)|
\ls \frac{2(n-1)}{n}\,\frac{nb_n}{n-1}\,S(Z)\,|Z|^{\frac{n-1}{n}}\end{equation}
by Lemma \ref{lem:zon-1}. Dividing by $|Z|^{\frac{n-1}{n}}$ we get the result. \prend

\medskip

\begin{question}\label{question-2}\rm It would be interesting to decide whether in the case of zonoids one has
\begin{equation}|Z|\,S(P_{\xi^{\perp }}(Z))\ls \frac{(n-1)}{n}\,S(Z)\,|P_{\xi^{\perp }}(Z)|\end{equation}
for every $\xi\in S^{n-1}$. This improvement of \eqref{eq:GHP} (for the class of zonoids) would give a sharp version of Theorem \ref{th:hyper-2}.
\end{question}

Next, we pass to estimates for the average surface area of hyperplane projections of $K$.

\begin{theorem}\label{th:hyper-3}Let $K$ be a convex body in ${\mathbb R}^n$. Then,
\begin{equation}\label{eq:hyper-9}|K|\,\int_{S^{n-1}}S(P_{\xi^{\perp }}(K))\,d\sigma (\xi )\ls \frac{2b_n}{n\omega_n^{\frac{1}{n}}}S(K)^2
\ls \frac{c_1}{\sqrt{n}}S(K)^2,\end{equation}
where $c_1>0$ is an absolute constant.
\end{theorem}

\noindent {\it Proof.} From \eqref{eq:GHP} we have
\begin{equation}\label{eq:hyper-10}|K|\,S(P_{\xi^{\perp }}(K))\ls \frac{2(n-1)}{n}\,S(K)\,|P_{\xi^{\perp }}(K)|\end{equation}
for every $\xi\in S^{n-1}$. Integrating on $S^{n-1}$ and using the identity
\begin{equation}\label{eq:hyper-11}S(K)=\frac{n\omega_n}{\omega_{n-1}}\int_{S^{n-1}}|P_{\xi^{\perp }}(K)|\,d\sigma (\xi )\end{equation}
we get
\begin{align}\label{eq:hyper-12}|K|\,\int_{S^{n-1}}S(P_{\xi^{\perp }}(K))\,d\sigma (\xi ) &\ls
\frac{2(n-1)}{n}\,S(K)\,\int_{S^{n-1}}|P_{\xi^{\perp }}(K)|\,d\sigma (\xi )\\
\nonumber &= \frac{2(n-1)}{n}\frac{\omega_{n-1}}{n\omega_n}S(K)^2\end{align}
Since
\begin{equation}\label{eq:hyper-13}\frac{2(n-1)}{n}\frac{\omega_{n-1}}{n\omega_n}=\frac{2b_n}{n\omega_n^{\frac{1}{n}}}\simeq \frac{1}{\sqrt{n}},\end{equation}
we get the result. \prend

\medskip

Now, let us assume that $K$ is in the minimal surface area, isotropic or John's position, or it is symmetric and in L\"{o}wner's position.
Then, from Proposition \ref{prop:surface-IJL} we know that
\begin{equation}\label{eq:hyper-14}S(K)\ls c_0n\,|K|^{\frac{n-1}{n}},\end{equation}
where $c_0>0$ is an absolute constant. From Theorem \ref{th:hyper-3} we get:

\begin{theorem}\label{th:hyper-4}Let $K$ be a convex body in ${\mathbb R}^n$. If $K$ is in the minimal surface area, isotropic or John's position, or it is symmetric and in L\"{o}wner's position, then
\begin{equation}\label{eq:hyper-15}|K|^{\frac{1}{n}}\,\int_{S^{n-1}}S(P_{\xi^{\perp }}(K))\,d\sigma (\xi )\ls c_2\sqrt{n}\,S(K)\end{equation}
where $c_2>0$ is an absolute constant.
\end{theorem}

\noindent {\it Note.} If $K$ is symmetric and in the minimal mean width position, using Proposition \ref{prop:surface-IJL} again, we get a weaker
(by a $\log n$ term) result:
\begin{equation}\label{eq:hyper-16}|K|^{\frac{1}{n}}\,\int_{S^{n-1}}S(P_{\xi^{\perp }}(K))\,d\sigma (\xi )\ls c_2\sqrt{n}(\log n)\,S(K)\end{equation}
where $c_2>0$ is an absolute constant.

\medskip

We pass now to lower bounds. Our analogue of \eqref{eq:intro-2} is the next theorem.

\begin{theorem}\label{th:hyper-5}Let $K$ be a convex body in ${\mathbb R}^n$. Then,
\begin{equation}\label{eq:hyper-17}\int_{S^{n-1}}S(P_{\xi^{\perp }}(K))\,d\sigma (\xi )\gr c_3\,S(K)^{\frac{n-2}{n-1}},\end{equation}
where $c_3>0$ is an absolute constant.
\end{theorem}

\noindent {\it Proof.} We write
\begin{align}\label{eq:hyper-18}\int_{S^{n-1}}S(P_{\xi^{\perp }}(K))\,d\sigma (\xi ) &=\frac{(n-1)\omega_{n-1}}{\omega_{n-2}}\int_{S^{n-1}}\int_{S^{n-1}\cap\xi^{\perp }}
|P_{\langle \xi ,\theta \rangle^{\perp }}(K)|\,d\sigma_{\xi^{\perp }}(\theta )\,d\sigma (\xi )\\
\nonumber &= \frac{(n-1)\omega_{n-1}}{\omega_{n-2}}\int_{G_{n,n-2}}|P_F(K)|\,d\nu_{n,n-2}(F).
\end{align}
From the Aleksandrov inequalities it follows that
\begin{align}\label{eq:hyper-19}\left (\frac{1}{\omega_{n-2}}\int_{G_{n,n-2}}|P_F(K)|\,d\nu_{n,n-2}(F)\right )^{\frac{1}{n-2}}
&\gr \left (\frac{1}{\omega_{n-1}}\int_{S^{n-1}}|P_{\xi^{\perp }}(K)|\,d\sigma (\xi )\right )^{\frac{1}{n-1}}\\
\nonumber &=  \left (\frac{S(K)}{n\omega_n}\right )^{\frac{1}{n-1}},\end{align}
which gives
\begin{equation}\label{eq:hyper-20}\int_{S^{n-1}}S(P_{\xi^{\perp }}(K))\,d\sigma (\xi )\gr
\frac{(n-1)\omega_{n-1}}{(n\omega_n)^{\frac{n-2}{n-1}}}\,S(K)^{\frac{n-2}{n-1}}\gr c_3S(K)^{\frac{n-2}{n-1}},\end{equation}
where $c_3>0$ is an absolute constant. \prend

\medskip

Now, let us assume that $K$ is in the minimal surface area, minimal mean width, isotropic, John or L\"{o}wner position.
Then, from Proposition \ref{prop:surface-IJL} (or from simple estimates in the cases of a not necessarily convex body $K$
that are not covered there) we know that, e.g.
\begin{equation}\label{eq:hyper-21}S(K)\ls c_0n^2\,|K|^{\frac{n-1}{n}},\end{equation}
where $c_0>0$ is an absolute constant. It follows that
\begin{equation}\label{eq:hyper-22}S(K)^{\frac{1}{n-1}}\ls (c_0n^2)^{\frac{1}{n-2}}\,|K|^{\frac{1}{n}}\ls c_4|K|^{\frac{1}{n}},\end{equation}
where $c_4>0$ is an absolute constant. Then,
\begin{equation}\label{eq:hyper-23}|K|^{\frac{1}{n}}\,S(K)^{\frac{n-2}{n-1}}\gr c_4^{-1}S(K).\end{equation}
Thus we have proved:

\begin{theorem}\label{th:hyper-6}Let $K$ be a convex body in ${\mathbb R}^n$. If $K$ is in the minimal surface area,
minimal mean width, isotropic, John's or L\"{o}wner's position, then
\begin{equation}\label{eq:hyper-24}|K|^{\frac{1}{n}}\,\int_{S^{n-1}}S(P_{\xi^{\perp }}(K))\,d\sigma (\xi )\gr c_5\,S(K),\end{equation}
where $c_5>0$ is an absolute constant.
\end{theorem}

\noindent {\it Note.} The proof of Theorem \ref{th:hyper-6} shows that \eqref{eq:hyper-24} continues to hold as long as
the mild condition
\begin{equation}S(K)^{\frac{1}{n-1}}\ls c|K|^{\frac{1}{n}}\end{equation}
is satisfied by $K$ with an absolute constant $c>0$.

\section{Quermassintegrals of hyperplane projections}

A generalization of \eqref{eq:GHP} to subspaces of arbitrary dimension and quermassintegrals of any order was given in \cite{Fradelizi-Giannopoulos-Meyer-2003}.

\begin{theorem}\label{th:FGM}Let $K$ be a convex body in ${\mathbb R}^n$ and let $0\ls p\ls k\ls n$. Then
for every $k$-dimensional subspace $F$ of ${\mathbb R}^n$, if $P_F(K)$ denotes the orthogonal projection of $K$ onto $F$, we have
\begin{equation}\label{eq:FGM}\frac{V_{n-p}(K)}{|K|}\gr\frac{1}{\binom{n-k+p}{n-k}}\frac{V_{k-p}(P_F(K))}{|P_F(K)|}.\end{equation}
\end{theorem}

Setting $k=n-1$, for every $1\ls p\ls n-2$ we have
\begin{equation}\label{eq:quer-1}\frac{V_{n-p}(K)}{|K|}\gr\frac{1}{p+1}\frac{V_{n-1-p}(P_{\xi^{\perp }}(K))}{|P_{\xi^{\perp }}(K)|}.\end{equation}
Therefore,
\begin{equation}\label{eq:quer-2}|K|\,\min_{\xi\in S^{n-1}}V_{n-1-p}(P_{\xi^{\perp }}(K))\ls (p+1)\,V_{n-p}(K)\,\min_{\xi\in S^{n-1}}|P_{\xi^{\perp }}(K)|,\end{equation}
and using \eqref{eq:rPiK} and Lemma \ref{lem:zon-1} we immediately get the following theorem.

\begin{theorem}\label{th:quer-1}Let $K$ be a convex body in ${\mathbb R}^n$. For every $1\ls p\ls n-2$ we have
\begin{equation}\label{eq:quer-3}|K|^{\frac{1}{n}}\,\min_{\xi\in S^{n-1}}V_{n-1-p}(P_{\xi^{\perp }}(K))\ls
\frac{(p+1)\omega_{n-1}\partial_K}{n\omega_n}\,V_{n-p}(K)\ls \frac{c_1(p+1)\partial_K}{\sqrt{n}}\,V_{n-p}(K),\end{equation}
where $c_1>0$ is an absolute constant. If $Z$ is a zonoid in ${\mathbb R}^n$ then, for every $1\ls p\ls n-2$ we have
\begin{equation}\label{eq:quer-4}|Z|^{\frac{1}{n}}\,\min_{\xi\in S^{n-1}}V_{n-1-p}(P_{\xi^{\perp }}(Z))\ls (p+1)b_n\,V_{n-p}(Z).\end{equation}
\end{theorem}

Starting from \eqref{eq:quer-1} and integrating on the sphere we get
\begin{align}\label{eq:quer-5}|K|\,\int_{S^{n-1}}V_{n-1-p}(P_{\xi^{\perp }}(K))\,d\sigma (\xi ) &\ls (p+1)V_{n-p}(K)\int_{S^{n-1}}|P_{\xi^{\perp }}(K)|\,d\sigma (\xi )\\
\nonumber &= \frac{(p+1)\omega_{n-1}}{n\omega_n}V_{n-p}(K)S(K).
\end{align}
Dividing by $|K|^{\frac{n-1}{n}}$ we get:

\begin{theorem}\label{th:quer-2}Let $K$ be a convex body in ${\mathbb R}^n$. For every $1\ls p\ls n-2$ we have
\begin{equation}\label{eq:quer-6}|K|^{\frac{1}{n}}\,\int_{S^{n-1}}V_{n-1-p}(P_{\xi^{\perp }}(K))\,d\sigma (\xi )\ls
\frac{(p+1)\omega_{n-1}}{n\omega_n}\frac{S(K)}{|K|^{\frac{n-1}{n}}}\,V_{n-p}(K).\end{equation}
In particular, if $K$ is in the minimal surface area, isotropic or John's position, or it is symmetric and in L\"{o}wner's position, then
we have
\begin{equation}\label{eq:quer-7}|K|^{\frac{1}{n}}\,\int_{S^{n-1}}V_{n-1-p}(P_{\xi^{\perp }}(K))\,d\sigma (\xi )\ls
c_1(p+1)\sqrt{n}\,V_{n-p}(K),\end{equation}
where $c_1>0$ is an absolute constant.
\end{theorem}

For the lower bound, an analogue of Theorem \ref{th:hyper-5}, we first observe that
\begin{equation}\label{eq:quer-8}V_{n-p}(K)=\omega_n\,[Q_{n-p}(K)]^{n-p}\quad\hbox{and}\quad
V_{n-1-p}(P_{\xi^{\perp }}(K))=\omega_{n-1}\,[Q_{n-1-p}(P_{\xi^{\perp }}(K))]^{n-1-p}\end{equation}
for every $\xi\in S^{n-1}$. Then, we write
\begin{align}\label{eq:quer-9}\int_{S^{n-1}}V_{n-1-p}(P_{\xi^{\perp }}(K))\,d\sigma (\xi ) &=\frac{\omega_{n-1}}{\omega_{n-1-p}}\int_{S^{n-1}}\int_{G_{\xi^{\perp }}(n-1,n-1-p)}|P_E(K)|\,d(E)\,d\sigma (\xi )\\
\nonumber &= \frac{\omega_{n-1}}{\omega_{n-1-p}}\int_{G_{n,n-1-p}}|P_F(K)|\,d\nu_{n,n-1-p}(F)\\
\nonumber &=\omega_{n-1}\,[Q_{n-1-p}(K)]^{n-1-p}.
\end{align}
From the Aleksandrov inequalities we have $Q_{n-1-p}(K)\gr Q_{n-p}(K)$, and hence
\begin{align}\label{eq:quer-10}\int_{S^{n-1}}V_{n-1-p}(P_{\xi^{\perp }}(K))\,d\sigma (\xi )&\gr \omega_{n-1}\,[Q_{n-p}(K)]^{n-2-p}\\
\nonumber &\gr \omega_{n-1}\left (\frac{V_{n-p}(K)}{\omega_n}\right )^{\frac{n-1-p}{n-p}}\end{align}
which gives the next theorem.

\begin{theorem}\label{th:quer-3}Let $K$ be a convex body in ${\mathbb R}^n$. For every $1\ls p\ls n-2$ we have
\begin{equation}\label{eq:quer-11}\int_{S^{n-1}}V_{n-1-p}(P_{\xi^{\perp }}(K))\,d\sigma (\xi )\gr
\frac{\omega_{n-1}}{\omega_n^{\frac{n-1-p}{n-p}}}\,[V_{n-p}(K)]^{\frac{n-1-p}{n-p}}.\end{equation}
\end{theorem}

Using the monotonicity of mixed volumes we may write
\begin{equation}\label{eq:quer-12}V_{n-p}(K)\ls V_{n-p}\big((K,n-p),(r(K)^{-1}K,p)\big)\ls \frac{|K|}{r(K)^p}.\end{equation}
Now, let us assume that $K$ is in the minimal surface area, isotropic or John's position, or it is symmetric and in L\"{o}wner's
position. Then, $r(K)\gr c_0|K|^{\frac{1}{n}}$ for an absolute constant $c_0>0$, and \eqref{eq:quer-12} gives
\begin{equation}\label{eq:quer-13}|K|^{\frac{n-p}{n}}\gr c_0^pV_{n-p}(K).\end{equation}
Therefore,
\begin{equation}\label{eq:quer-14}|K|^{\frac{1}{n}}[V_{n-p}(K)]^{\frac{n-1-p}{n-p}}\gr c_0^{\frac{p}{n-p}}[V_{n-p}(K)]^{\frac{1}{n-p}}[V_{n-p}(K)]^{\frac{n-1-p}{n-p}}
=c_0^{\frac{p}{n-p}}\,V_{n-p}(K).\end{equation}
From Theorem \ref{th:quer-3} we get:

\begin{theorem}\label{th:quer-4}Let $K$ be a convex body in ${\mathbb R}^n$, which is in the minimal surface area, isotropic or John's position,
or it is symmetric and in L\"{o}wner's position. For every $1\ls p\ls n-2$ we have
\begin{equation}\label{eq:quer-15}|K|^{\frac{1}{n}}\int_{S^{n-1}}V_{n-1-p}(P_{\xi^{\perp }}(K))\,d\sigma (\xi )\gr
\frac{\omega_{n-1}c_0^{\frac{p}{n-p}}}{\omega_n^{\frac{n-1-p}{n-p}}}\,V_{n-p}(K)\gr \left (\frac{c_1}{n}\right )^{\frac{p}{2(n-p)}}\,V_{n-p}(K),\end{equation}
where $c_1>0$ is an absolute constant.
\end{theorem}

Note that $\left (\frac{c_1}{n}\right )^{\frac{p}{2(n-p)}}\gr c_2$ for an absolute constant $c_2>0$ as long as $p\ls cn/(\log n)$.

\section{Surface area of projections of higher codimension}

Recall that $mV_{m-1}(A)=S(A)$ for every convex body $A$ in ${\mathbb R}^m$. Therefore, setting $p=1$ in Theorem \ref{th:FGM} we get:

\begin{lemma}\label{lem:higher-1}Let $K$ be a convex body in ${\mathbb R}^n$ and let $1\ls k\ls n-1$. Then
for every $k$-dimensional subspace $F$ of ${\mathbb R}^n$ we have
\begin{equation}\label{eq:higher-1}\frac{S(K)}{|K|}\gr \frac{n}{k(n-k+1)}\frac{S(P_F(K))}{|P_F(K)|}.\end{equation}
\end{lemma}

We first prove an analogue of Theorem \ref{th:hyper-2}.

\begin{theorem}\label{th:higher-2}Let $Z$ be a zonoid in ${\mathbb R}^n$ and let $1\ls k\ls n-1$. Then,
\begin{equation}\label{eq:higher-2}|Z|^{\frac{n-k}{n}}\,\min_{F\in G_{n,k}}S(P_F(Z))\ls
(n-k+1)b_n^{n-k}\,S(Z).\end{equation}
\end{theorem}

\noindent {\it Proof.} From \eqref{eq:higher-1} we see that
\begin{align}\label{eq:higher-3}|Z|\,\min_{F\in G_{n,k}}S(P_F(Z)) &\ls \frac{k(n-k+1)}{n}\,S(Z)\,\min_{F\in G_{n,k}}|P_F(Z)|
\ls \frac{k(n-k+1)}{n}\,\frac{nb_n^{n-k}}{k}\,S(Z)\,|Z|^{\frac{k}{n}}\\
\nonumber &= (n-k+1)b_n^{n-k}S(Z)\,|Z|^{\frac{k}{n}},\end{align}
where in the last step we have also used Theorem \ref{th:zon-2}. Dividing by $|Z|^{\frac{k}{n}}$ we get the result. \prend

\begin{definition}\label{def:higher-4}\rm For every convex body $K$ in ${\mathbb R}^n$ and every $1\ls k\ls n-1$ we introduce
the parameter
\begin{equation}\label{eq:higher-4}p_k(K):=\frac{1}{|K|^{\frac{k}{n}}}\int_{G_{n,k}}|P_F(K)|\,d\nu_{n,k}(F).\end{equation}
\end{definition}

Using Lemma \ref{lem:higher-1} and applying the same argument as in the proof of Theorem \ref{th:hyper-3} we get:

\begin{theorem}\label{th:higher-5}Let $K$ be a convex body in ${\mathbb R}^n$ and let $1\ls k\ls n-1$. Then,
\begin{equation}\label{eq:higher-5}|K|^{\frac{n-k}{n}}\,\int_{G_{n,k}}S(P_E(K))\,d\nu_{n,k}(E)\ls \frac{k(n-k+1)}{n}\,S(K)\,p_k(K).\end{equation}
\end{theorem}

\noindent {\it Proof.} From Lemma \ref{lem:higher-1} we have
\begin{equation}\label{eq:higher-6}|K|\,S(P_F(K))\ls \frac{k(n-k+1)}{n}\,S(K)\,|P_F(K)|\end{equation}
for every $F\in G_{n,k}$. Integrating with respect to $\nu_{n,k}$ on $G_{n,k}$ we get
\begin{equation}\label{eq:higher-7}|K|\,\int_{G_{n,k}}S(P_F(K))\,d\nu_{n,k}(F)\ls \frac{k(n-k+1)}{n}\,S(K)\,\int_{G_{n,k}}|P_F(K)|\,d\nu_{n,k}(F),\end{equation}
and the result follows. \prend

\begin{remark}\rm Let us assume that $K$ is in the minimal surface area, isotropic or John's position, or it is symmetric and in L\"{o}wner's
position. From \eqref{eq:quer-12} and the fact that $r(K)\gr c_0|K|^{\frac{1}{n}}$ we get
\begin{equation}p_k(K)=\frac{1}{|K|^{\frac{k}{n}}}\frac{\omega_k}{\omega_n}V_k(K)\ls \frac{1}{|K|^{\frac{k}{n}}}\frac{\omega_k}{\omega_n}\frac{|K|}{c_0^{n-k}|K|^{\frac{n-k}{n}}}
=\frac{\omega_k}{\omega_nc_0^{n-k}}.\end{equation}
Then, Theorem \ref{th:higher-5} gives the following analogue of Theorem \ref{th:hyper-4}:
\end{remark}

\begin{theorem}Let $K$ be a convex body in ${\mathbb R}^n$. If $K$ is in the minimal surface area, isotropic or John's position, or it is symmetric and in L\"{o}wner's position, then
\begin{equation}|K|^{\frac{n-k}{n}}\,\int_{G_{n,k}}S(P_F(K))\,d\nu_{n,k}(F)\ls \frac{k(n-k+1)}{n}\frac{\omega_k}{\omega_nc_0^{n-k}}\,S(K)\end{equation}
where $c_0>0$ is an absolute constant.
\end{theorem}

The lower bound of Theorem \ref{th:hyper-5} can be generalized as follows.

\begin{theorem}\label{th:higher-6}Let $K$ be a convex body in ${\mathbb R}^n$. For every $1\ls k\ls n-1$ we have
\begin{equation}\label{eq:higher-8}\int_{G_{n,k}}S(P_F(K))\,d\nu_{n,k}(F)\gr \frac{k\omega_k}{(n\omega_n)^{\frac{k-1}{n-1}}}\,S(K)^{\frac{k-1}{n-1}}.\end{equation}
\end{theorem}

\noindent {\it Proof.} We write
\begin{align}\label{eq:higher-9}\int_{G_{n,k}}S(P_F(K))\,d\nu_{n,k}(F) &=\frac{k\omega_k}{\omega_{k-1}}\int_{G_{n,k}}\int_{E\cap\xi^{\perp }}
|P_{F\cap \xi^{\perp }}(K)|\,d\sigma_{F}(\xi )\,d\nu_{n,k} (F)\\
\nonumber &= \frac{k\omega_k}{\omega_{k-1}}\int_{G_{n,k-1}}|P_E(K)|\,d\nu_{n,k-1}(E).
\end{align}
From the Aleksandrov inequalities we have
\begin{align}\label{eq:higher-10}\left (\frac{1}{\omega_{k-1}}\int_{G_{n,k-1}}|P_E(K)|\,d\nu_{n,k-1}(E)\right )^{\frac{1}{k-1}}
&\gr \left (\frac{1}{\omega_{n-1}}\int_{S^{n-1}}|P_{\xi^{\perp }}(K)|\,d\sigma (\xi )\right )^{\frac{1}{n-1}}\\
\nonumber &=  \left (\frac{S(K)}{n\omega_n}\right )^{\frac{1}{n-1}},\end{align}
which gives
\begin{equation}\label{eq:higher-11}\int_{G_{n,k}}S(P_F(K))\,d\nu_{n,k}(F)\gr
\frac{k\omega_k}{(n\omega_n)^{\frac{k-1}{n-1}}}\,S(K)^{\frac{k-1}{n-1}},\end{equation}
as claimed. \prend

\medskip

Now, let us assume that $K$ is in the minimal surface area, isotropic or John's position, or it is symmetric and in L\"{o}wner's position.
Then, from Proposition \ref{prop:surface-IJL} we know that
\begin{equation}\label{eq:higher-12}S(K)\ls c_0n\,|K|^{\frac{n-1}{n}},\end{equation}
where $c_0>0$ is an absolute constant. Therefore,
\begin{equation}\label{eq:higher-13}|K|^{\frac{n-k}{n}}S(K)^{\frac{k-1}{n-1}}\gr \frac{1}{(c_0n)^{\frac{n-k}{n-1}}}S(K)^{\frac{n-k}{n-1}}S(K)^{\frac{k-1}{n-1}}
=\frac{1}{(c_0n)^{\frac{n-k}{n-1}}}S(K),\end{equation}
and Theorem \ref{th:higher-6} implies the following.

\begin{theorem}\label{th:higher-7}Let $K$ be a convex body in ${\mathbb R}^n$. If $K$ is in the minimal surface area, isotropic or John's position, or it is symmetric and in L\"{o}wner's position, then
\begin{equation}\label{eq:higher-14}|K|^{\frac{n-k}{n}}\,\int_{G_{n,k}}S(P_F(K))\,d\nu_{n,k}(F)\gr
\frac{k\omega_k}{(n\omega_n)^{\frac{k-1}{n-1}}(c_0n)^{\frac{n-k}{n-1}}}\,S(K),\end{equation}
where $c_0>0$ is an absolute constant.
\end{theorem}

\bigskip

\bigskip

\noindent {\bf Acknowledgements.} The second named author was partially supported by the US National Science Foundation grant DMS-1265155.

\bigskip

\bigskip

\footnotesize
\bibliographystyle{amsplain}

\bigskip

\bigskip

\noindent \textsc{Apostolos \ Giannopoulos}: Department of
Mathematics, University of Athens, Panepistimioupolis 157-84,
Athens, Greece.

\smallskip

\noindent \textit{E-mail:} \texttt{apgiannop@math.uoa.gr}

\bigskip

\noindent \textsc{Alexander \ Koldobsky}: Department of
Mathematics, University of Missouri, Columbia, MO 65211.

\smallskip

\noindent \textit{E-mail:} \texttt{koldobskiya@missouri.edu}

\bigskip

\noindent \textsc{Petros \ Valettas}: Department of
Mathematics, University of Missouri, Columbia, MO 65211.

\smallskip

\noindent \textit{E-mail:} \texttt{valettasp@missouri.edu}

\end{document}